\documentclass{article}
\usepackage[parfill]{parskip}    		
\usepackage{amssymb}
\usepackage{amsmath}
\usepackage{amsthm}
\usepackage{url}
\usepackage{natbib}
\usepackage[pdftex]{graphicx}

\usepackage{booktabs}
\usepackage{url}

\title{Robustness of semiparametric efficiency in \\ nearly-true models for two-phase samples}
\author{Thomas Lumley}

\begin{document}
\maketitle
\begin{abstract}
  We examine the performance of efficient and AIPW estimators under two-phase sampling when the complete-data model is nearly correctly specified, in the sense that the misspecification is not reliably detectable from the data by any possible diagnostic or test.  Asymptotic results for these nearly true models are obtained by representing them as sequences of misspecified models that are mutually contiguous with a correctly specified model.  We find that for the least-favourable direction of model misspecification the bias in the efficient estimator induced can be comparable to the extra variability in the AIPW estimator, so that the mean squared error of the efficient estimator is no longer lower. This can happen when the most-powerful test for the model misspecification still has modest power.  We verify that the theoretical results agree with simulation in three examples: a simple informative-sampling model for a Normal mean, logistic regression in the classical case-control design, and linear regression in a two-phase design\\
\end{abstract}

\section{Introduction}
For parametric models there is a substantial literature on when estimates or predictions from a slightly-misspecified submodel are more accurate than those from a correctly-specified full model. Considerations of bias--variance tradeoff in using a deliberately misspecified model date back at least to the question of linear vs quadratic discriminant analysis, and when it was worth using linear discriminant analysis even though predictor variances were probably different in the classes being discriminated.  Simulations showed that linear discriminant analysis tended to be preferable for small samples and higher dimensions, with quadratic discrimination winning out as sample size increased \citep{gilbert-lda, marks-dunn,wahl-kronmal}.   \citet{hjort-tness} was perhaps the first to give precise analytical results for a question of this sort.  \citet[Chapter 5]{claeskens-hjort} review and expand on this literature.

In this paper we consider a similar question in a particular class of semiparametric models for incomplete data.   When are design-based estimates from a model assuming only known sampling probabilities more accurate than the semiparametric-efficient estimator that also assumes a specific parametric or semiparametric regression model? It is important to note that we do not address the issue of data missing by happenstance, where the sampling probabilities cannot be assumed known.  We consider only the problem of two-phase sampling where measurements can be taken reliably on the selected subsample: for example, by assays on previously frozen blood samples, by expert coding of existing free-text responses, or by paying for data from commercial databases.

The case of gross model misspecification is relatively uninteresting, since a competent analyst motivated by asymptotic efficiency would then choose a different model. Heuristically, since the difference between naive and efficient estimators is typically $O_p(n^{-1/2})$, the interesting bias:variance tradeoffs happen at biases of that order. We would like to consider \emph{nearly true} models, defined as those where tests of model misspecification will not reliably reject.  

In applied statistics this would refer to some practical set of tests of model specification.
For asymptotic purposes we use a stricter definition of `nearly true', that even the Neyman--Pearson test comparing the sequence of data-generating distributions to the closest sequence of distributions in the model has power bounded away from one. That is, since we do not know which forms of misspecification are most common in practice, we will consider worst-case misspecification.  We are able to give simple quantitative results for the impact of worst-case misspecification. In simple but realistic examples, worst-case misspecification of the model wipes out the advantage of the semiparametric-efficient estimator when most powerful goodness-of-fit test has power less than 75\% (and as low as 25\%) at 5\% size.

\citet{watson-holmes} and their discussants review recent Bayesian decision theory research based on a local-minimax approach similar to ours. They assume a prior and likelihood are specified, but then consider decisions that are minimax over all distributions within a given (small) Kullback--Leibler divergence from the posterior.

We begin by outlining the estimators and models we are considering, and the definition of the target of inference when the model is misspecified.  Assuming that efficient design-based and model-based estimators are available, we then derive the asymptotic distribution of these estimators under contiguous model misspecification. We present three examples where the asymptotic results can be verified in simulation, and then discuss the interpretation of the results. Code for all the simulations is available from \url{https://github.com/tslumley/nearlytrue}.  Finally, we discuss the implications of the results and their sensitivity to choices made in the setup.

\section{Incomplete data models}
\label{notation}
Our estimators are based around an outcome model ${\cal P}$ for $[Y|X,Z]$, where $Y$ and $Z$  are available for everyone, but $X$ may be available only for a subsample of individuals.  This model is indexed by the parameter of interest $\theta$ and a typically infinite-dimensional nuisance parameter $\eta$.  

With complete data on $(X, Y, Z)$, estimation of $\theta$ would be performed by solving
$$U(\theta)=U(\theta; X, Y,Z; \hat\eta)=0,$$ 
where $U(\theta)$ is an estimate of the efficient score or efficient influence function for $\theta$ giving at least locally efficient estimation with complete data at the true $\eta_0$. We call the resulting estimator $\tilde\theta$ the complete-data estimator; in survey statistics it is known as the census parameter.
We {\em  define} the target of inference $\theta^*$ as the large-sample limit of the complete-data estimator, and we use this definition whether or not the true distribution of $[Y|X,Z]$ satisfies the outcome model.

The complete-data estimator is not available in practice, because the data are incomplete. 
We obtain $X$ only for a sample of the observations. We also observe other variables $L$ that are not part of the efficient influence function for $\theta$ with complete data, and refer to $(Z, Y, L)$ as the `phase 1 data '. We write $R_i$ for the indicator that $X$ is observed for observation $i$, and assume that $\pi_i=E[R_i|\textrm{phase 1 data}]$ is known (this can be weakened to  `known up to an identifiable finite-dimensional parameter')

We define two semiparametric models for the observed data
\begin{enumerate}
\item The \emph{sampling-only model} ${\cal S}$ is the model for $(X,Y,Z,L,R)$ that assumes only 
$$E[R_i|X_i,Y_i,L_i,Z_i]=\pi_i$$
\item The \emph{regression submodel} ${\cal M}$ is the submodel where the conditional distribution $[Y|X,Z]$ satisfies the outcome model, so ${\cal M} = {\cal S} \cap {\cal P}$
\end{enumerate}

We describe an estimator that is consistent for $\theta^*$ under the sampling-only model as a design-based estimator, and a estimator that is consistent for $\theta^*$ under the regression submodel but not the sampling-only model as a model-based estimator.

 A simple design-based estimator is the Horvitz--Thompson or Inverse-Probability Weighted (IPW) estimator, which solves
\begin{equation}
\sum_{i=1}^N \frac{R_i}{\pi_i} U(\theta; X_i, Y_i, Z_i,\hat\eta)=0.
\label{eq:ht}
\end{equation}

The Horvitz--Thompson estimator does not use any information from $(Z,Y,L)$ on observations with $R_i=0$. \citet{rrz} defined Augmented IPW estimators (AIPW) that solve
\begin{equation}
\sum_{i=1}^N \frac{R_i}{\pi_i} U(\theta; X, Y, Z,\hat\eta) + \sum_{i=1}^N \left(1-\frac{R_i}{\pi_i}\right)A(\theta; Z_i,Y_i,L_i,\hat\phi)=0
\label{eq:aipw}
\end{equation}
where $A$ is an arbitrary function of phase-one data.  The most efficient choice of $A_i(\theta)$ is 
$$A_i^*=E[U_i(\theta)|\textrm{phase 1 data}],$$
 though this will often not be feasible.  In practice, reasonably efficient choices of $A_i$ are conveniently available via connections with calibration of weights in survey sampling \citep{han-mi-calib,breslow-SiB-twophase,breslow-AJE-twophase,samuelsen-stoer,lumley-calibration,robins-cal-aipw,deville-raking}.

RRZ also characterized the efficient estimator, which uses a complete-data influence function $V_i$ obtained by projecting $U$ orthogonal to the tangent space for nuisance parameters
\begin{equation}
\sum_{i=1}^N \frac{R_i}{\pi_i} V(\theta; X_i, Y_i, Z_i) + \sum_{i=1}^N \left(1-\frac{R_i}{\pi_i}\right)A^\dag_i(\theta)=0\label{eq:eff}
\end{equation}
with $$A^\dag_i(\theta)=E[V_i(\theta)|\textrm{phase 1 data}].$$

If the distribution is in $\cal P$ and $V$ can be estimated sufficiently accurately then $\hat\theta_{\mathrm{eff}}$ is (locally) semiparametric efficient, and typically is strictly more efficient than the best AIPW estimator.
 We will write $\check{U}$ and $\check{V}$ for the influence functions of the two estimators  
 $$\check{U}(\theta)=\frac{R_i}{\pi_i} U(\theta; X_i, Y_i, Z_i)+\left(1-\frac{R_i}{\pi_i}\right)A_i^*(\theta; Z_i,Y_i,\hat\phi)$$ 
 and 
$$\check{V}(\theta)= \frac{R_i}{\pi_i} V(\theta; X_i, Y_i, Z_i) +  \left(1-\frac{R_i}{\pi_i}\right)A_i^\dag(\theta )=0.$$

The characterization of the efficient estimator given by RRZ is not necessarily the most convenient way to compute the efficient estimator. For example, computations using profile likelihoods are described by \citet{scott-wild-06} for the estimators proposed by those authors and co-workers.  For our purposes it is sufficient that equation~\ref{eq:eff} characterizes the semiparametric-efficient estimator up to asymptotic equivalence.  We do not need to assume that equation~\ref{eq:eff} is used in implementation, and (except in performing simulations) we do not need to know the efficient influence function $V$ explicitly. 

There are many important technical issues in constructing an efficient estimator that we will not cover in this paper since we are assuming that such an estimator has in fact been constructed.  Modern discussions of many of these issues can be found in \citet{tsiatis-book} and \citet{kosorok-book}.

\section{Nearly-true models}
\label{definition}
\subsection{Definition}

The practical question for data analysis underlying the concept of a nearly-true model is whether it is sufficient to conduct tests or examine diagnostics for model misspecification in order to justify relying on the efficient estimator.   This question leads to the heuristic concept of a nearly-true model as a model that cannot reliably be rejected by the available diagnostics.  Since our tools for proving statements about efficiency are asymptotic, we need a formal characterization of  `nearly true' that captures this heuristic concept but allows relevant asymptotic arguments to be constructed.

The available tools for model criticism will vary by the model and the data collected, but a bound on the effectiveness of these tools is given by the Neyman--Pearson lemma.  If we knew that the data came either from a specific distribution $P_{\theta,\nu}$ inside the model or a specific distribution $Q$ outside the model, the most powerful test would be based on the likelihood ratio $L=dQ/dP_{\theta,\nu}$.     We can thus measure the distance from $Q$ to the model using $\inf_{\theta,\nu} dQ/dP_{\theta,\nu}$, the minimum Kullback--Leibler divergence.

For any fixed $Q$ and $P_{\theta,\nu}$ the test based on $L$ will eventually reject with certainty as $N$ and $n$ increase.  To construct an asymptotic setting that is relevant to the practical question we need a sequence $Q_n$ of misspecified distributions where $dQ_n/dP_{\theta,\nu,n}$ is bounded.   That is, the data at hand are considered as an element of a sequence of experiments in which $Q_n$  is not reliably distinguishable from the model.  

A formal characterisation of this condition is that the sequence of data distributions $Q_n$ and some sequence of model distributions $P_n$ are mutually contiguous [eg Chapter 6, van der Vaart, 1998]. The definition of mutual contiguity is that for any sequence of events $A_n$, 
$$Q_n[A_n]\to 0 \iff P_n[A_n]\to 0.$$
In particular, this holds if $A_n$ is the event that we find a satisfactory level of model fit in the sample of size $n$ after using some set of diagnostics. 

When $Q_n$ and $P_n$ are mutually contiguous the sequence of likelihood ratios $L_n=dQ_n/dP_n$
is uniformly tight. If  this sequence  converges in distribution under $P_n$ to a variable $L_\infty$ then $E[L_\infty]=1$.  By taking a subsequence if necessary it is no loss of generality to assume that this convergence in distribution holds.

\subsection{Estimation in nearly true models}
\label{lecam}

LeCam's Third Lemma \citep{lecam-contiguity,vdv} describes how to relate distributions of sequences of statistics $D_n$ under models $P_n$ to their distributions under $Q_n$.  In this section we will write $D_n\stackrel{P_n}{\rightsquigarrow}D$ to mean $D_n\stackrel{d}{\to}D$ under $P_n$ and $D_n\stackrel{Q_n}{\rightsquigarrow}D$ to mean $D_n\stackrel{d}{\to}D$ under $Q_n$.

We are concerned only with the case when $D_n$ and the log likelihood $\log L_n = \log dP_n/dQ_n$ are asymptotically multivariate Normal. 

\paragraph{Lemma}
If
$$\left(D_n,\, \log\frac{dQ_n}{dP_n}\right)\stackrel{P_n}{\rightsquigarrow} N\left( \left(\begin{array}{c}\mu\\-\kappa^2/2\end{array}\right),\;\left(\begin{array}{cc} \Sigma & \tau\\ \tau^T & \kappa^2\end{array}\right)\right)$$
then
$$D_n\stackrel{Q_n}{\rightsquigarrow}N(\mu+\tau,\,\Sigma).$$

That is, the change from $P_n$ to $Q_n$ shifts the limiting distribution but does not change the scale.

 In particular, if $\mu=0$ and $D$ is scalar we can write $\sigma^2$ for $\Sigma$ and reparametrize $\tau$ in terms of a correlation $\tau=\rho\kappa\sigma$.   
 We then have
$$D_n\stackrel{P_n}{\rightsquigarrow} N(0, \sigma^2)$$
and 
$$D_n\stackrel{Q_n}{\rightsquigarrow} N(\kappa\rho\sigma,\sigma^2)$$

Here $\rho$ is the correlation between $\log L_\infty$ and $D$ under $P_n$. It describes whether the model is misspecified in a direction that affects $\theta$. The size of the model misspecification, in terms of the power of the most powerful test for misspecification,  is measured by $\kappa$. In order to construct the misspecification giving the worst bias, we need to reduce $\theta$ to a univariate parameter of primary interest. If $\theta$ is multivariate, this may require taking a single element or a linear combination of elements.  From now on we assume $\theta$ is univariate.

 If we take 
$$D_n=\sqrt{n}(\hat\theta_{\mathrm{eff}}-\hat\theta_{AIPW})$$
and
$$\sqrt{n}(\hat\theta_{\mathrm{eff}}-\hat\theta_{AIPW})\stackrel{P_n}{\rightsquigarrow} N(0,\omega^2)$$
LeCam's third lemma gives
$$\sqrt{n}(\hat\theta_{\mathrm{eff}}-\hat\theta_{AIPW})\stackrel{Q_n}{\rightsquigarrow} N(\kappa\rho\omega,\omega^2)$$

Under $Q_n$ the outcome model is misspecified, so care is needed in defining the `true' parameter value.  We define $\theta^*$ as the value to which the outcome-model point estimator would converge with complete data as $N\to\infty$. 

The AIPW estimator is still asymptotically unbiased, so
$$\sqrt{n}(\hat\theta_{\mathrm{AIPW}}-\theta^*)\stackrel{Q_n}{\rightsquigarrow} N(0,\sigma^2+\omega^2)$$
but
$$\sqrt{n}(\hat\theta_{\mathrm{eff}}-\theta^*)\stackrel{Q_n}{\rightsquigarrow} N(\kappa\rho\omega,\sigma^2)$$

Finally, we note that the likelihood ratio test for $H_0:\;Q_n$ vs $H_1:\,P_n$ prescribed by the Neyman--Pearson Lemma has null distribution $N(-\kappa^2/2,\kappa^2)$ and alternative distribution $N(\kappa^2/2,\kappa^2)$, so it is equivalent to a (one-sided) test for a location shift of $\kappa$ in a $N(0,1)$ distribution.  Its power at level 0.05 is 13\% for $\kappa=0.5$, 26\% for $\kappa=1$, 64\% for $\kappa=2$, and 90\% for $\kappa=3$.  For $\kappa\leq 3$ the Neyman--Pearson test could certainly not be described as reliable, and in most scenarios the available model diagnostics will be less powerful than the Neyman--Pearson test: the alternative will not be known precisely, and tests based on adding parameters to the model will typically be two-sided, if not multivariate.

When $\rho\neq 0$, $\hat\theta_{\mathrm{eff}}$ is asymptotically biased for $\theta^*$. If $\hat\theta_{\mathrm{AIPW}}$ is locally efficient among AIPW estimators there will exist sequences $Q_n$ with $\rho$ arbitrarily close to 1.  For AIPW estimators other than the most efficient one $\rho$ will be bounded away from one; in particular, $\rho$ will typically be bounded away from 1 for the Horvitz--Thompson estimator.

The asymptotic mean squared error of $\hat\theta_{\mathrm{eff}}$ is 
$$MSE_{\mathrm{eff}}=\kappa^2\rho^2\omega^2+\sigma^2$$
and of  $\hat\theta_{\mathrm{AIPW}}$ is
$$MSE_{\mathrm{AIPW}}=\sigma^2+\omega^2$$
If $\kappa^2\rho^2>1$, $\hat\theta_{\mathrm{AIPW}}$ has smaller mean squared error.

For the best AIPW estimator there are misspecified models $Q_n$ with  $\rho$ arbitrarily close to 1, so small amounts of model misspecification in an unfavourable direction are sufficient to remove the advantage of the efficient estimator.  For the crude Horvitz--Thompson estimator, on the other hand, the maximum attainable $\rho$ may be quite small and the efficient estimator may have substantially superior mean-squared error when the model is nearly correct.

\section{Simulation study}
\label{simul}
A simulation example to verify the theoretical results presented above requires the ability to compute the efficient estimator, an AIPW estimator that is close to optimal, and the efficient influence function; and the ability to sample from a misspecified distribution with $\rho\gg 0$. We present three such examples: estimation of a simple Normal mean from subsampled data, logistic regression in the classical nested case--control design, and linear regression under two-phase outcome-dependent sampling.

We consider both simple but \emph{ad hoc} examples of misspecification, and attempts to approximate the worst-case distribution by exponential tilting.  That is, given $\check{V}$ and $\check{U}$, and a density $f_0$ in the model, we construct
$$f_\kappa\propto f_0\exp\left[\kappa(\check{V}-\check{U})\right]$$
an exponential family through $f_0$. Writing $C_\kappa$ for the normalising constant, the score function for $\kappa$ is
$$\frac{\partial\log f}{\partial \kappa} = \check{V}-\check{U} + \frac{\partial \log C_\kappa}{\partial\kappa}.$$ 
To the extent that the last term is small, this family accurately approximates the worst-case model misspecification. 

These examples also allow us to examine the reasonableness of the worst-case misspecifications. This is important: the direction to `generic' misspecifications in a high- or infinite-dimensional space may be nearly orthogonal to the worst-case direction. If the worst-case misspecification appears too pathological, its practical relevance will be diminished.

\subsection{Normal mean}
\label{normal-mean}
We start with a toy example of informative sampling with no phase-1 information. 
Let $X_i\sim N(\mu,1)$, for $i=1,2,\dots,n$, and suppose $(X_i,\pi_i)$ is observed with probability 
$$\pi_i = \frac{e^{x_i}}{1+e^{x_i}}$$
Let $R_i$ be the indicator that observation $i$ is observed, and let
$$p_0(\mu)=\int  \frac{1}{1+e^{x_i}}\phi(x-\mu)\,dy$$
be the marginal probability of not being observed, where $\phi(\cdot)$ is the standard Normal density. 

The log likelihood is 
$$\log L(\mu) = \sum_{i=1}^n R_i\log\phi(x-\mu) + (1-R_i)\log p_0(\mu)$$
The maximum likelihood estimator is efficient in this parametric model. It is not available in closed form, but is easy to compute numerically.  The (efficient) score function is 
$$\check V_i(x;\mu) = R_i(x_i-\mu)+(1-R_i) p_0'(\mu)$$

The efficient AIPW estimator is simply the Horvitz--Thompson estimator, since there is no auxiliary complete-data information to augment it with, and the weighted estimating function is
$$\check U(x;\mu) = \frac{R_i}{\pi_i}(x_i-\mu)$$
The target of inference under the sampling-only model, $\theta^*$, is the population mean of $X$. 

\begin{figure}
\includegraphics[width=\textwidth]{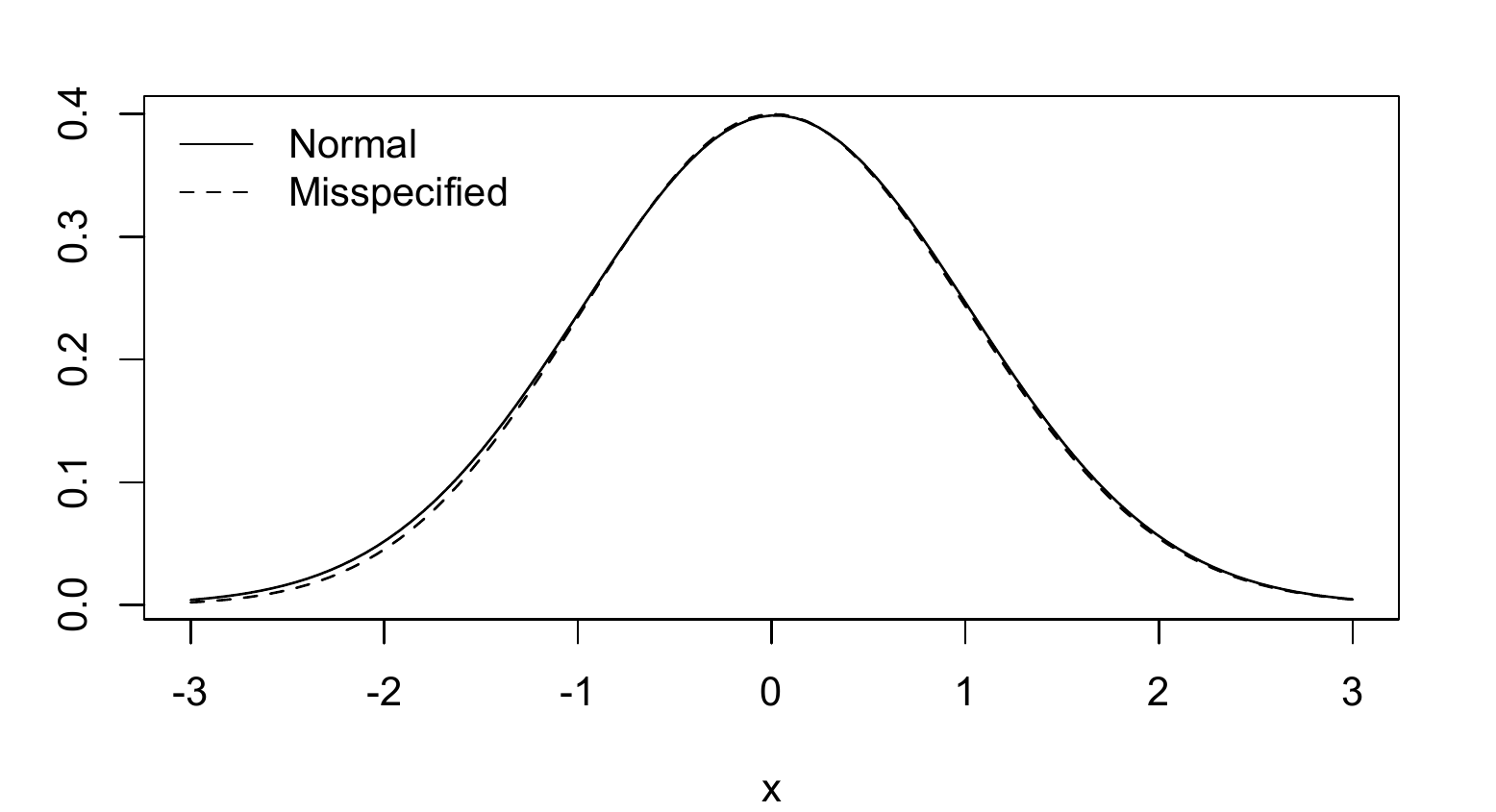}
\caption{Densities of assumed Normal model and worst-case misspecified model with equal mean squared error for estimating the population mean, $n=10^4,\,\delta=0.012$.}
\label{toy-density}
\end{figure}

The worst-case misspecification is approximated by 
$$f_{\tilde Q_\kappa}(x;\mu) = C_\kappa \phi(x-\mu)\exp\left[ -\kappa(2-e^{-x})(x-\mu)-\kappa p_0'(\mu) \right]$$
with the approximation improving as $\kappa\to 0$.  Here $C_\kappa$ is a normalising constant that can be computed by numerical integration, and $\tilde Q_\kappa$ can be simulated by rejection sampling from a Normal distribution. $\tilde Q_\kappa$ has a slightly higher peak than a Normal distribution, and asymmetrically lighter shoulders, as shown in Figure~\ref{toy-density}.  Note in particular that, in contrast to familiar examples of poor outlier-resistance of the Normal MLE, the misspecified model is not heavy-tailed.

\begin{figure}
\includegraphics[width=\textwidth]{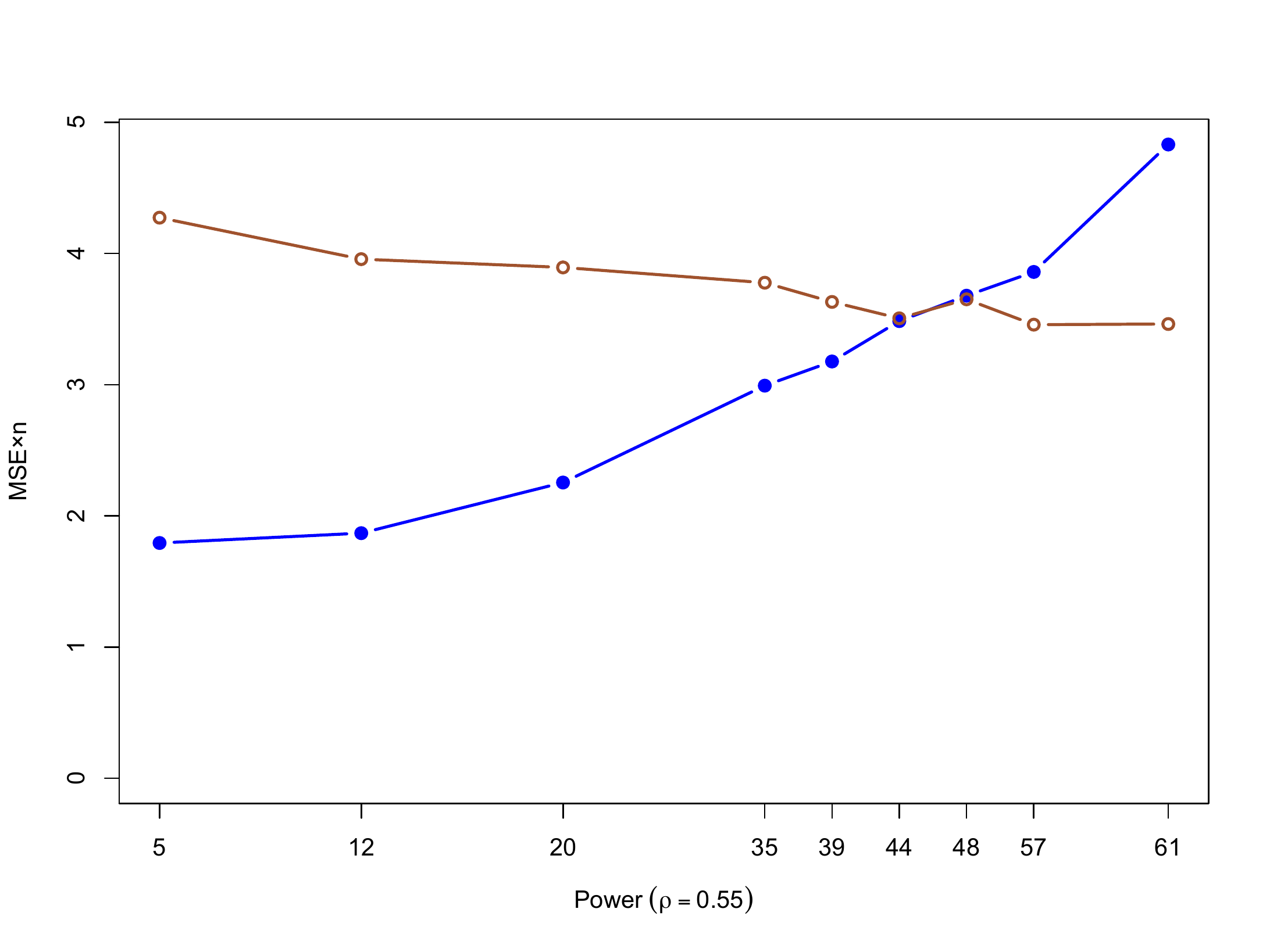}
\caption{Mean squared error for weighted estimator (open circles) and misspecified MLE (filled circles) in the Normal mean subsampling model, against power of the Neyman--Pearson test. $n=10^4,\,\delta=0.012$, $\rho\approx 0.55$, $10^4$ simulations for each set of values.}
\label{toy-sim}
\end{figure}

At $\mu=0$ and with the model correctly specified, the Horvitz--Thompson estimator has 37\% efficiency, i.e., $\omega^2\approx 1.7\sigma^2$. The asymptotic construction of the worst-case misspecification gives $\rho \approx 0.55$ at $n=10^4$.  Figure~\ref{toy-sim} shows the mean squared error for the efficient and weighted estimators.  When the mean squared error is equal for the two estimators, the power of the misspecification test is about 50\%, so the misspecification is not reliably detectable. Equal mean squared error occurs at $\kappa\approx 0.5$ and $\kappa\rho\approx 0.9$, in reasonable agreement with the asymptotic analysis.


\subsection{Case-control design}
 
An important incomplete-data design for which all the required quantities are known is the classical population-based case--control design.  In phase one a binary outcome variable $Y$ is measured on a large sample, and in phase two, predictors $X$ are measured on all case subjects with $Y=1$ and on a fraction $\pi_0$ of control subjects with $Y=0$.   We consider the case of a single predictor $X$.
The model is 
$$E[Y|X=x]=\mu(\theta,x)=\frac{e^{\alpha+x\beta}}{1+e^{\alpha+x\beta}}$$

The complete data efficient influence function is
$$U(\theta)= E[XX^T\mu(\theta,X)(1-\mu(\theta,X)] X\frac{1}{\pi}(Y-\mu(\theta,X))$$
so the IPW estimator solves
$$\sum_{i=1}^n \frac{1}{\pi_i}X_i(Y_i-\mu(\theta, X))=0$$
and
$$\check{U}(\theta)=E[XX^T\frac{\mu(\theta,X)(1-\mu(\theta,X)}{\pi}|R=1] X(Y-\mu(\theta,X))$$
and because there is no further phase-one information this is also the best AIPW estimator, ie, $A^*(\theta)=0$.

The efficient estimator for $\beta$ is the unweighted case--control estimator \citep{prentice-pyke,breslow-wellner-robins-2000}, with influence function
$$\check{V}(\beta)= E[XX^T\tilde\mu(\theta,X)(1-\tilde\mu(\theta,X)|R=1] X(Y-\tilde\mu(\theta,X))$$
where 
$$\tilde\mu=E[Y|X=x,R=1]=\frac{e^{\alpha-\log\pi_0+x\beta}}{1+e^{\alpha-\log\pi_0+x\beta}}$$ is the regression function conditional on being sampled in phase two.

\begin{figure}
\centerline{\includegraphics[height=3in]{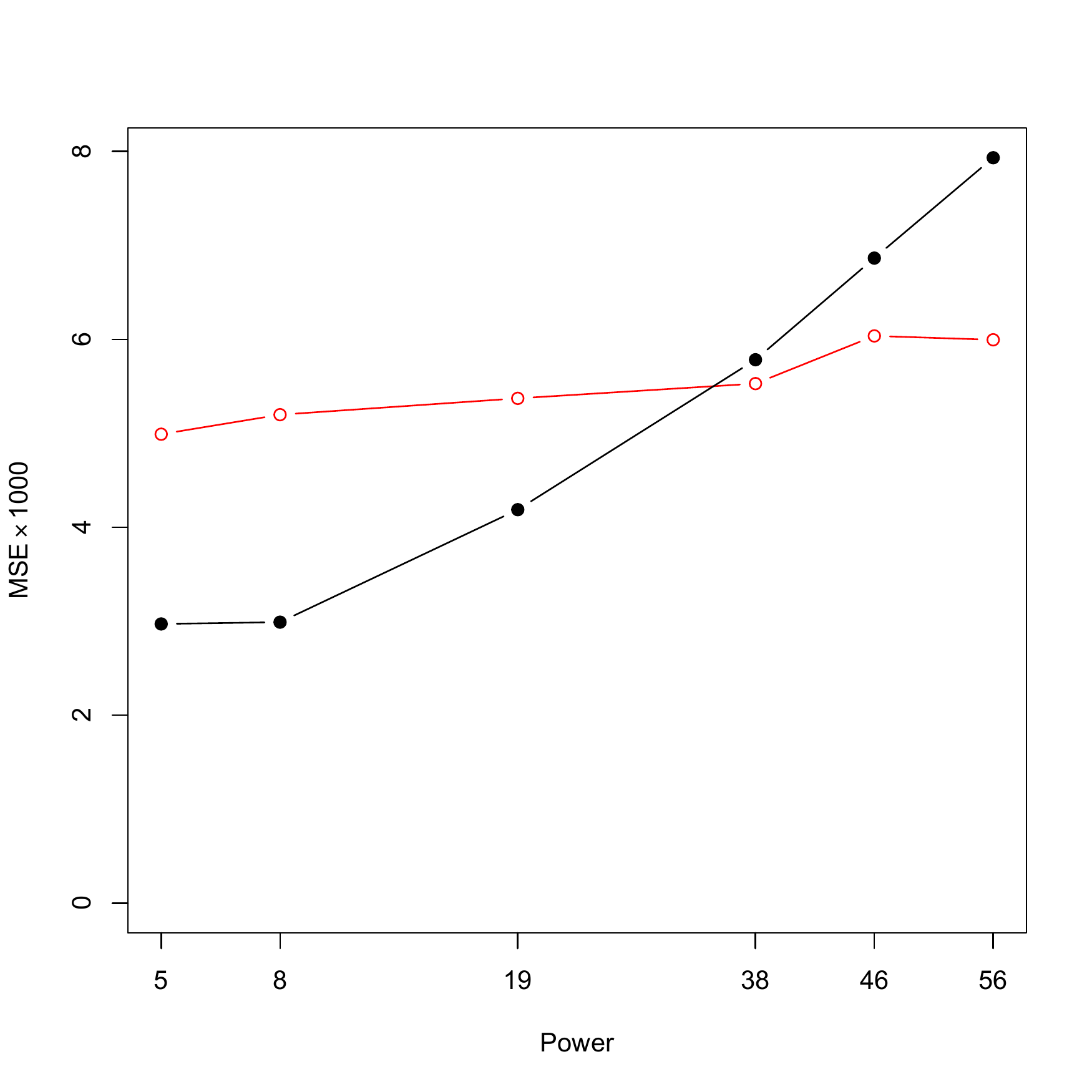}}

\caption{Efficiency of MLE (filled circles) and design-based case--control estimator (open circles) of $\beta$ under quadratic misspecfication $(\rho\approx 0.75)$, with power for the efficient test}
\label{cc-quadratic}
\end{figure}

\citet{scott-wild-02} compared weighted and unweighted estimation with case--control sampling under a misspecified model: 
$$\mathrm{logit}\,\Pr[Y=1]=\alpha-\log\pi_0+\beta X +\gamma X^2$$
 where $X\sim N(0,1)$ and $\gamma$ was chosen to give power in the vicinity of 50\% for rejecting $\gamma=0$ with the standard Wald test.  As a starting point, we consider the same model. It is worth noting that the loss of efficiency and susceptibility to bias are larger for Normal $X$ than for many other covariate distributions: the qualitative picture will be the same in other settings, but the difference between the weighted and unweighted estimators will be smaller.  Figure~\ref{cc-quadratic} shows the mean squared error for the two estimators in terms of the power of the one-sided test for the quadratic term.   For $\beta=1$, the quadratic has $\rho\approx 0.75$.   .

As a second example, we approximated the worst-case model by resampling with exponential tilting.  That is, first we generated a population of size $2\times 10^5$ with $X\sim N(0,1)$ and $\mathrm{logit}\,\Pr(Y=1|X=x)=-5+X$.  We then took all cases and an equal number of controls for an intermediate sample; let $m$ be the number of cases.  On the intermediate sample we computed $\Delta(x,y)=\check{V}(x,y)-\check{U}(x,y)$.  We then sampled $m/2$ cases from the cases in the intermediate sample and $m/2$ controls from the controls in the intermediate sample, with replacement, with sampling probabilities for individual $i$ proportional to $\exp(\epsilon\Delta(x_i,y_i)$. We repeated this procedure $10^4$ times for each value of $\epsilon$.
 
The likelihood ratio conditional on the realisation of the finite population is the product of the resampling probabilities for the selected observations.  We verified that the log likelihood ratios were well approximated by $N(\kappa^2/2,\kappa^2)$ with $\kappa$ depending on $\epsilon$, as the asymptotic arguments show.  The correlation between the log likelihood ratio and $\hat\beta-\tilde\beta$ was $\rho\approx 0.5$. 

Since the likelihood ratio test is also sensitive to the distribution of $X$, which would not typically be of interest to the analyst, we conducted a test for misspecification of the conditional distributions of $Y|X$.  A data set 100 times larger than the individual simulated data sets was constructed by concatenating 100 realisations, and both a straight line and a smoothing spline with 4 degreees of freedom were fitted for the relationship between $\mathrm{logit}\,\Pr(Y=1|X=x)$ and $X$ using the {\sf gam} package\citep{package-gam}.  A new simulated data set of the original size was constructed and the Bernoulli loglikelihood computed for the fitted spline curve and fitted straight line from the large data set; the difference in these two loglikelihoods gives the Neyman--Pearson test statistic of these two alternatives for the conditional distribution.  The conditional test statistic was then computed 1000 times. 

\begin{figure}
\centerline{\includegraphics[height=3in]{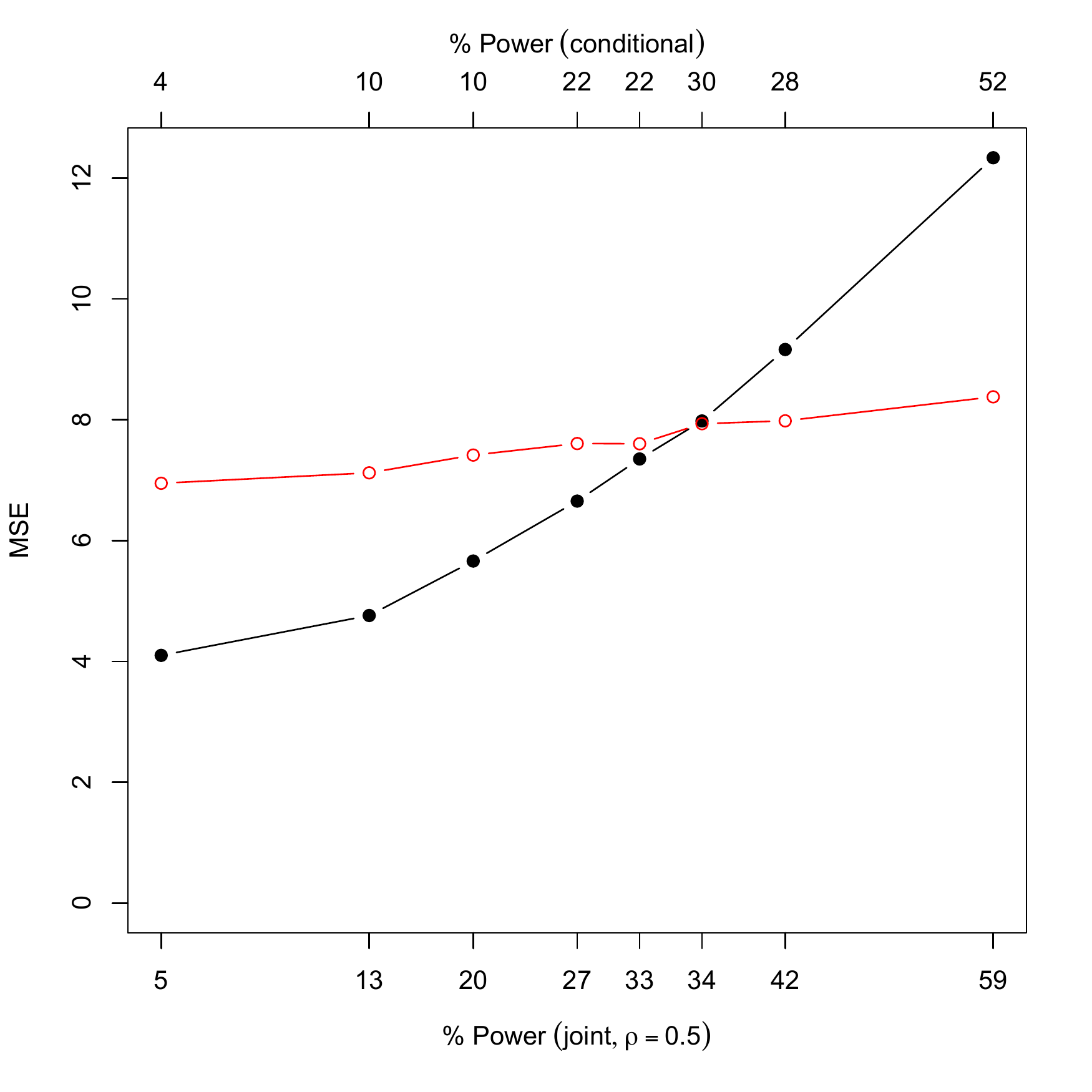}}
\caption{Efficiency of MLE (filled circles) and design-based case--control estimator (open circles) of $\beta$ with misspecification by exponential tilting. Lower $x$-axis shows power of test for misspecification in the joint model for $(X,Y)$, with $\rho\approx 0.5$, upper $x$-axis for misspecification in the model for $Y|X$.}
\label{cc-efficiency}
\end{figure}

Figure~\ref{cc-efficiency} shows the mean squared error for the weighted estimator and the efficient estimator as a function of power. The lower $x$-axis shows the power of the test based on the resampling probabilities; the upper $x$-axis shows the power for the test of the conditional distribution.  The misspecification is closer to the true worst case for the conditional distributions: the power of the test at $\kappa=1$ is only 30\%.  The joint test, because it is also sensitive to the difference in the distribution of $X$, is more powerful, but not in a way that is useful in typical data analysis because it is not usual to rely on assumptions about the distribution of the predictors in logistic regression.

\begin{figure}[p]
\centerline{\includegraphics[height=3in]{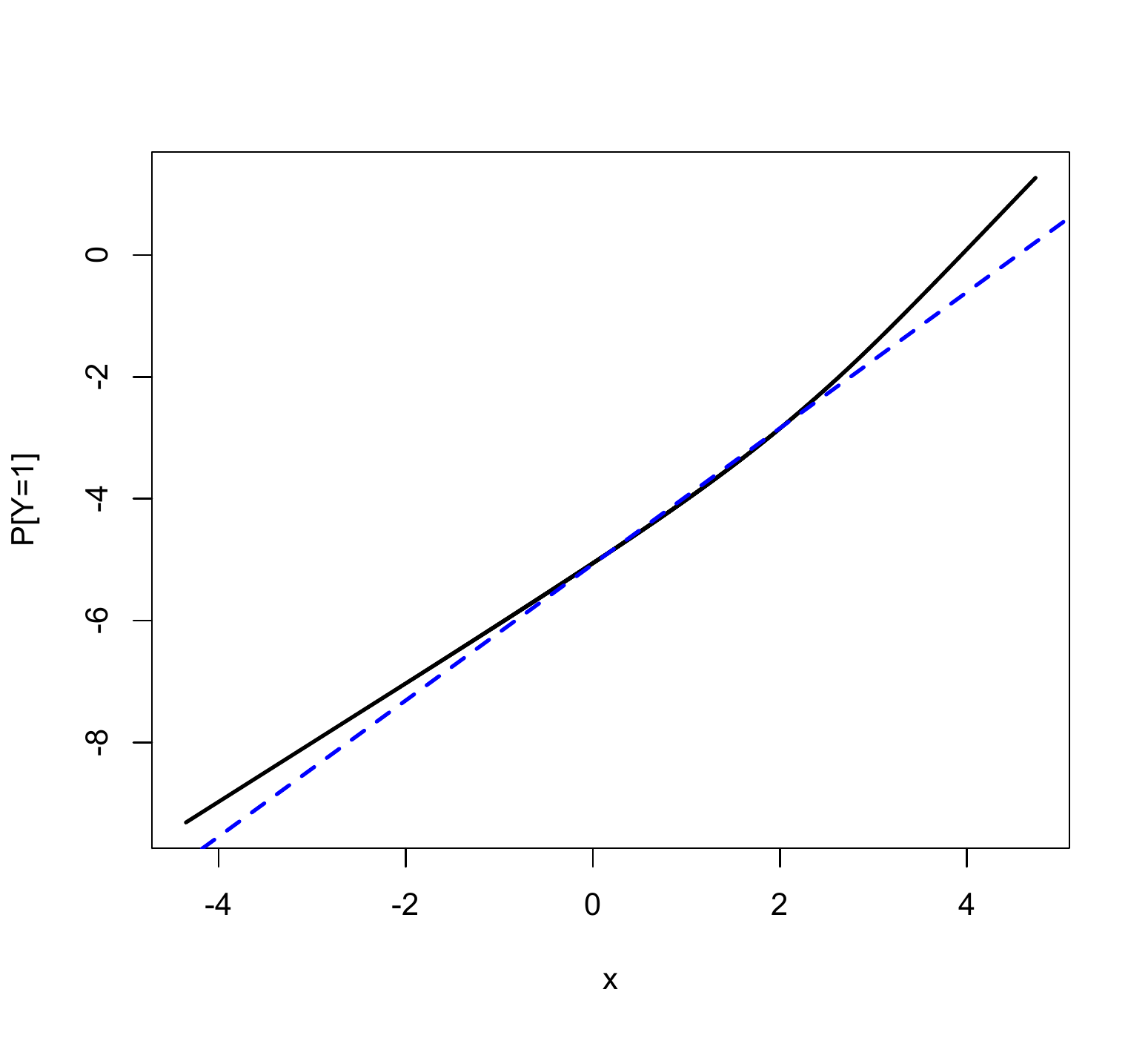}}
\caption{Approximately worst-case misspecified model for $\mathrm{logit}\Pr[Y=1|X=x]$ with $\kappa\approx 1$, 2000 cases, $X\sim N(0,1)$ and best-fitting straight line.}
\label{cc-distortion}
\end{figure}

Figure~\ref{cc-distortion} shows the fitted smoothing spline and fitted straight line for $\mathrm{logit}\Pr[Y=1|X=x]$, based on the large simulated sample. The deviation from linearity is not quadratic; the curvature is greater near the mean of $X$ and the function becomes approximately linear further out.

\begin{figure}[p]
\centerline{\includegraphics[height=3in]{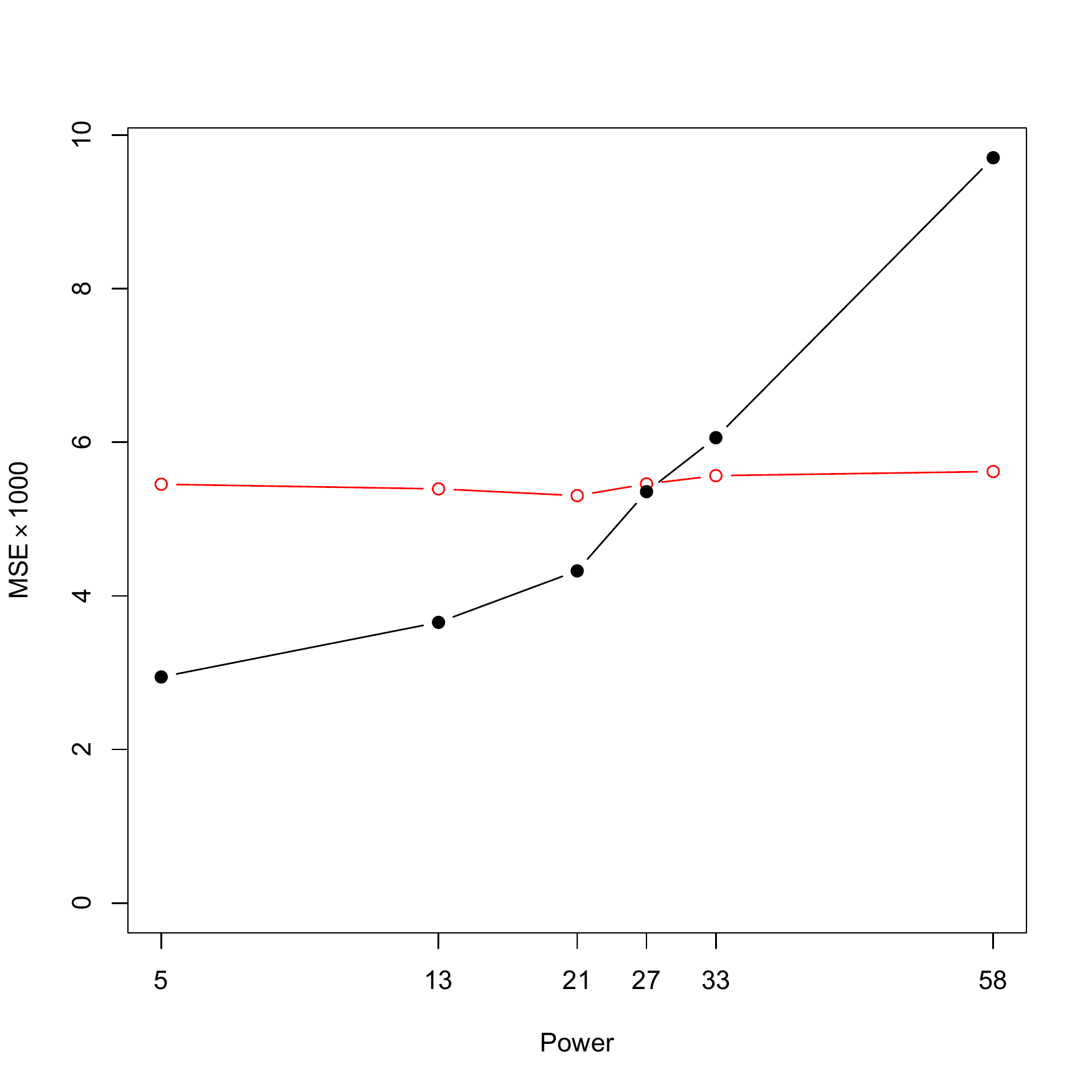}}
\caption{Efficiency of MLE (filled circles) and design-based case--control estimator (open circles) of $\beta$ under linear-spline misspecfication $(\rho\approx 0.92)$ with a single knot at $x=1.2$. The $x$-axis shows power for the efficient test}
\label{cc-linspline}
\end{figure}

Finally, motivated by Figure~\ref{cc-distortion}, we considered a linear-spline misspecification. First, we chose the location of the single knot at $x=1.8$ by maximising the correlation between the misspecification test and the bias of the misspecified mle.  Then, in independent simulations, we varied the size of the misspecification.  The linear spline (Figure~\ref{cc-linspline}) is very close to the worst-case misspecification, having $\rho=0.92$.

\subsection{Linear regression under two-phase sampling}

Our third example is of linear regression under two-phase sampling.   We consider a simple model with
$$Y=X+e$$
where $X$ and $\epsilon$ are standard Normal. In addition, we have an auxiliary variable $Z=X+e^*$, with $e^*$ also standard Normal, observed for everyone.

At the first phase of sampling, we observe the outcome variable $Y$ and the auxiliary variable $Z$ for an iid sample of 4000 individuals. At the second phase we measure $X$ for all individuals with $Z>2.33$ and $Z<2.33$ (the expected 5th and 95th percentiles) and for 200 with $Z\in [-2.33,\,2.33]$.  The rationale for this form of sampling is to increase efficiency by oversampling the influential points; an additional benefit when this is done in medical research is that individuals extreme on one variable are often worth oversampling for their values of other variables. 

We want to fit the model
$$E[Y|X=x]=\alpha+\beta x$$
The semiparametric maximum likelihood estimator of $(\alpha,\beta)$ when $Y$ is Normal is available in the {\sf missreg3} package for R\citep{package-missreg3}, using the methods described by\citep{scott-wild-06}.  Following~\citep{breslow-AJE-twophase}, we construct a reasonably efficient design-based estimator by 
\begin{enumerate}
\item Fitting an imputation model to predict $X$ from $Z$ and $Y$ for all observations 
\item Fitting $E[Y|\hat X=\hat x] = \tilde\alpha+\tilde\beta\hat x$ with the imputed data
\item Using the influence functions from the the model as auxiliary variables to calibrate the sampling weights.
\item Fitting $E[Y|X=x]=\alpha+x\beta$ to the phase-2 subsample using the calibrated weights.
\end{enumerate}

The {\sf missreg3} package does not directly provide the influence functions, so we approximated them by the empirical influence functions, the changes in the estimates when a single observation is deleted.  Exploratory analysis suggested that the misspecification should be approximately linear over the region $Z\in[-2.33,2.33]$ where the data are incomplete. 

\begin{figure}
\centerline{\includegraphics[height=3in]{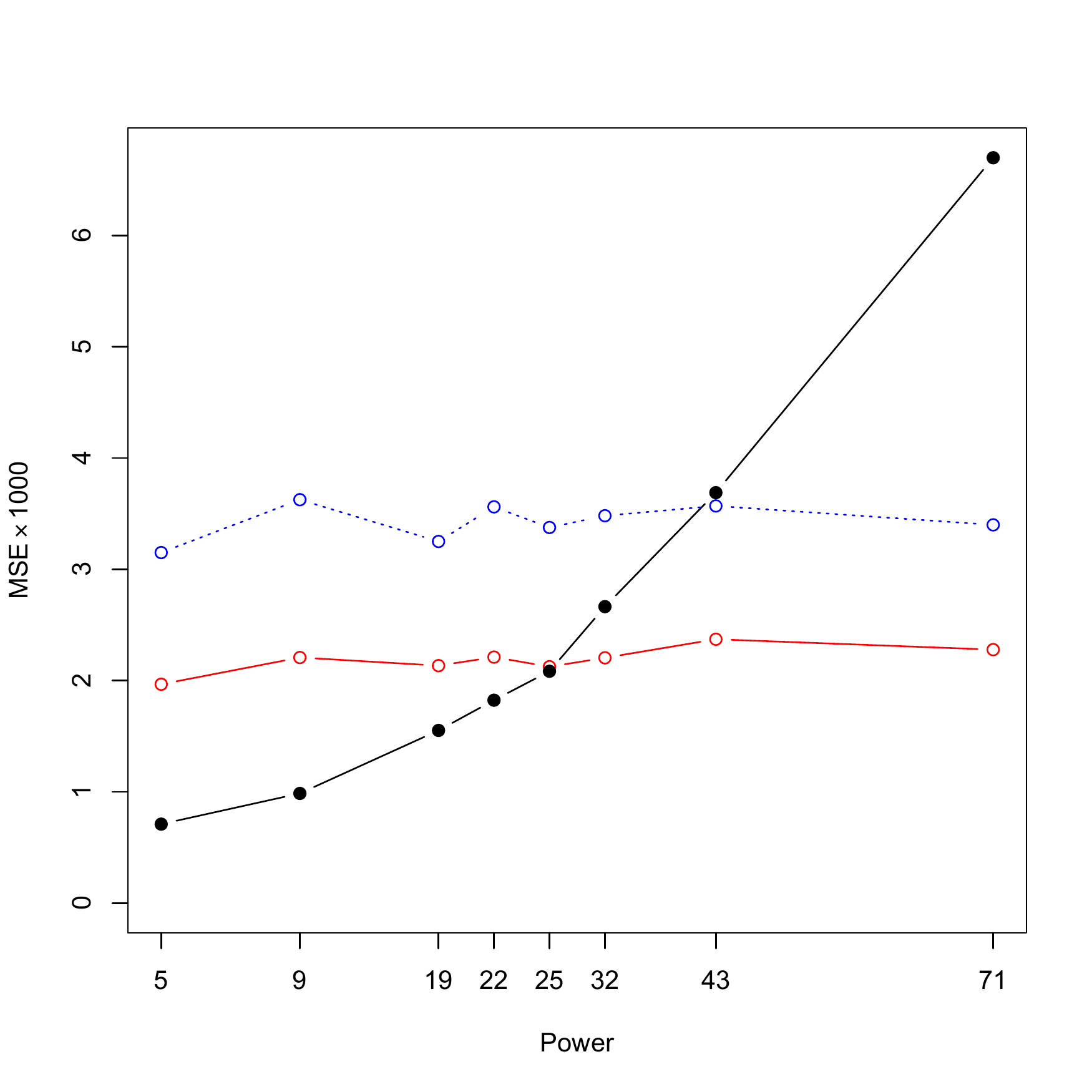}}

\caption{Efficiency of MLE (filled circles), calibrated design-based estimator (open circles), and Horvitz--Thompson estimator (dashed lines) of $\beta$ for the linear regression model under two-phase sampling. Model misspecfication is linear in $x$ over the incompletely-sampled region of $Z$, and has $\rho\approx 0.7$ . The $x$-axis shows power for the efficient test of misspecification}
\label{linear}
\end{figure}

Simulation using a linear model misspecification over this region,
$$E[Y|X=x]=\alpha+\beta x+\gamma x\times \{z \in [-2.33,2.33]\},$$
gives $\rho\approx 0.7$, and~Figure~\ref{linear} shows the mean squared errors from the efficient estimator, the  calibrated design-based estimator, and the Horvitz--Thompson estimator. Again, it is possible to construct misspecifications where the design-based estimator has lower mean squared error but the test for model misspecification (the one-sided test of $\gamma=0$) is not reliable. 

This example differs from the previous two in that the Horvitz--Thompson estimator is not the efficient design-based estimator, and its performance is noticeably worse than the calibrated estimator.

\section{Discussion}
\label{discussion}
In a specific setting of regression models under two-phase sampling, we have shown that even an unrealistically powerful goodness-of-fit test cannot reliably diagnose model misspecification that is serious enough to remove the efficiency advantage of a semiparametric-efficient estimator.    Our result is in some ways a converse of \citet{freedman-diagnostics}. His paper shows that diagnostics can only be powerful against narrowly-specified alternatives.  We see here that even against narrowly-specified alternatives the power need not be sufficient to validate important assumptions that are not already known to be true.  We do not advocate the routine use of pre-testing, which can seriously degrade the distribution of an estimator \citep{guggenberger-twostage,guggenberger-panel,leeb-poetscher}, but the low power of the Neyman--Pearson test is still useful as a indicator of closeness.

The results depend critically on two choices: the definition of the target parameter as the large-sample limit of the complete-data estimator, and the use of worst-case contiguous alternatives. The first choice is standard in survey statistics and a traditional option in biostatistics. As \citet{mantel-haenszel-1959} observed, \emph{A primary goal is to reach the same conclusions in a retrospective study as would have been obtained from a forward study, if one had been done}.  The second choice fits with the standard definition of semiparametric efficiency, which uses the same local asymptotic minimax characterisation of efficiency within the model that we use outside the model.

Our results are entirely asymptotic; perhaps more importantly, our simulations have relied on large sample sizes.  We expect the same phenomena to exist at more usual sample sizes, but sampling from a model close to the worst-case becomes more computationally difficult when $\delta n^{-1/2}$ is larger. 

An early reviewer of the material in this paper claimed that our conclusions are obvious in light of the Local Asymptotic Minimax Theorem. In a heuristic sense this is arguably true, but our quantitative results do not follow straightforwardly from the Local Asymptotic Minimax Theorem. The theorem, \citep[section 8.8]{vdv} gives a bound for the error in an arbitrary sequence of estimators $T_n$ of a parameter $\zeta$.  If $\ell()$ is any `bowl-shaped' loss function and the efficient estimator $\hat\zeta$ satisfies $\sqrt{n}(\hat\zeta-\zeta_0)\stackrel{d}{\to}\mathbb{Z}$ then
\begin{equation}
\lim_{\delta\to\infty}\liminf_{n\to\infty} \sup_{\zeta: \sqrt{n}|\zeta-\zeta_0|\leq\delta} E_\zeta\left[\ell(\sqrt{n}(T_n-\zeta))\right]\geq E\ell(\mathbb{Z})
\label{lam}
\end{equation}

Now let $\zeta_0$ be a point in the regression model, but define $\mathbb{Z}$ in terms of the efficient estimator in the sampling-only model. An efficient estimator in the regression model is just another estimator $T_n$ of $\zeta$, and so cannot beat the efficient sampling-only estimator uniformly over $\delta n^{-1/2}$-neighbourhoods of $\zeta_0$, at least as $\delta\to\infty$.   Our result is stronger, though narrower: if $T_n$ is efficient under the regression model we give an explicit, finite value for $\delta$ that makes the inequality in display \ref{lam} strict, for a loss that picks out a particular one-dimensional contrast.

In the setting of case--control sampling, \citet{scott-wild-02} argue that the efficient estimator may be useful even though it is biased. This is a reasonable point of view: the case--control study is a special situation both because the bias is analytically relatively tractable and because there is so much practical experience with the design.   There may also be other situations where it makes sense to use the efficient estimator even when it is biased.  One scenario is when the primary interest is in testing rather than  estimation and the bias does not affect the null hypothesis. For example, in a case--control design, it is possible to test the null hypothesis that $Y$ is independent of $X$ using a likelihood ratio test, because if $Y$ and $X$ are independent the logistic regression model with $\beta=0$ will be correctly specified.   This specific example is not compelling because the design-based estimator is fully efficient when $\beta=0$, so there is no increase in power for large samples and contiguous alternatives. If there is a difference in power in small samples it would need to be demonstrated directly and would not follow automatically from the greater efficiency of the MLE.   The efficiency bounds for the Cox model under case--cohort sampling \citep[Figure 3]{nan-info-bound} suggest that AIPW estimators have full or nearly full efficiency at the null hypothesis in this setting as well, though \citet{claudia-calibration} suggests this is not the case for Cox regression under countermatching.

There may also be situations where there are good reasons to believe the key assumptions of the outcome model are sufficiently close to being true. Independence assumptions justified by Mendelian segregation in genetics or by randomisation in clinical trials, and linear or power relationships derived from physics are two examples.  Such an argument must rely on substantive knowledge of a particular application, rather than on the observed data.  Conversely, there are many situations where the sampling probabilities that we have assumed to be known may be misspecified because of non-response.  Misspecified sampling probabilities would lead to bias in the AIPW estimator, and the conservative modelling strategy would be to use a doubly-robust estimator \citep{bang-robins,kang-doubly-robust,han-mi-calib}

In many two-phase designs we do not have either analytic results or sufficient experience to trust intuition in handling misspecification bias.   The results for contiguous models show that correct model specification is effectively an unverifiable assumption at the level of precision at which  discussions of relative efficiency take place. More study of this phenomenon is needed both to characterise the behaviour under realistic kinds of misspecification and in order to understand when it is appropriate to accept the changed target of estimation in order to increase precision.

\citet{mukherjee-chatterjee}, in a related problem, describe a shrinkage strategy between efficient and robust estimation. While tempting, this strategy fails to be adaptive precisely under local misspecification of the sort we consider, at least when the two estimators have the same convergence rate; when they have different convergence rates adaptation is possible \citep{mehdi-rates}.

There is sometimes a substantial difference in efficiency between the crude Horvitz--Thompson estimator and computationally simple AIPW estimators based on calibration of weights, although the gain will be small if the available auxiliary data are not very predictive.  In contrast to the use of the semiparametric efficient estimator, using an improved AIPW estimator, at least in large enough samples, is a `free lunch': there is a gain in precision with no change in assumptions.  This fact, combined with our results on misspecification, suggest that
simulation studies of efficient estimators under two-phase sampling should at least consider contiguous model misspecification and should compare to a more efficient AIPW estimator rather than the  Horvitz--Thompson estimator where possible.

\bibliographystyle{abbrvnat}
\bibliography{../survey,../raoscott}

\end{document}